\documentclass[12pt,letter]{amsart}
\usepackage{amsmath,amssymb,amsthm,amscd,amsxtra,amsfonts,mathrsfs,graphicx,enumerate,bm,slashed,array,mathtools}
\input xy
\xyoption{all}
\newtheorem{Theorem}{Theorem}[section]
\newtheorem{Lemma}[Theorem]{Lemma} 
\newtheorem{Proposition}[Theorem]{Proposition}
\newtheorem{Corollary}[Theorem]{Corollary}
\newtheorem{Example}[Theorem]{Example}
\newtheorem{Remark}[Theorem]{Remark}
\newtheorem{Question}[Theorem]{Question}
\newtheorem{Definition}[Theorem]{Definition}

\newtheorem{Physics Remark}[Theorem]{Physics Remark}

\newtheorem*{Theorem A}{Theorem A}
\newtheorem*{Theorem C}{Theorem C}
\newtheorem*{Definition B}{Definition B}
\newtheorem*{Definition D}{Definition D}
\newtheorem*{Theorem D}{Theorem D}
\newtheorem*{Theorem E}{Theorem E}
\newtheorem*{Theorem F}{Theorem F}
\advance\evensidemargin-.5in
\advance\oddsidemargin-.5in
\advance\textwidth1in

\begin{document}
 \author{Charlie Beil}
 \address{Heilbronn Institute for Mathematical Research, School of Mathematics, Howard House, The University of Bristol, Bristol, BS8 1SN, United Kingdom.}
 \email{charlie.beil@bristol.ac.uk}
 \title{Nonnoetherian geometry}
 \keywords{Non-noetherian rings, foundations of algebraic geometry, noncommutative algebraic geometry.}
 \subjclass[2010]{13C15, 14A20}
 \date{}

\begin{abstract}
We introduce a theory of geometry for nonnoetherian commutative algebras with finite Krull dimension.
In particular, we establish new notions of normalization and height: \textit{depiction} (a special noetherian overring) and \textit{geometric codimension}.
The resulting geometries are algebraic varieties with positive dimensional points, and are thus inherently nonlocal.
These notions also give rise to new equivalent characterizations of noetherianity that are primarily geometric. 

We then consider an application to quiver algebras whose simple modules of maximal dimension are one dimensional at each vertex.
We show that the vertex corner rings of $A$ are all isomorphic if and only if $A$ is noetherian, if and only if the center $Z$ of $A$ is noetherian, if and only if $A$ is a finitely generated $Z$-module. 
Furthermore, we show that $Z$ is depicted by a commutative algebra generated by the cycles in its quiver.
We conclude with an example of a quiver algebra where projective dimension and geometric codimension, rather than height, coincide.
\end{abstract}

\maketitle

\tableofcontents

\section{Introduction}

The purpose of this paper is to introduce a new framework for understanding the geometry of nonnoetherian algebras with finite Krull dimension.
The underlying idea we propose is that nonnoetherian geometry should be the geometry of \textit{nonlocal} algebraic varieties.

We first motivate our study of nonnoetherian geometry in the following applications.

(1) \textit{Non-cancellative dimer algebras.}  
A dimer algebra is a type of quiver algebra whose quiver embeds in a torus, with homotopy-like relations.
These noncommutative algebras encode certain gauge theories in string theory.
Dimer algebras with the cancellation property are Calabi-Yau algebras, and their centers are 3-dimensional Gorenstein singularities.  
Non-cancellative dimer algebras, on the other hand, are not Calabi-Yau, and their centers are nonnoetherian.  
In contrast to their cancellative counterparts, very little is known about these algebras, despite the fact that almost all dimer algebras are non-cancellative.
In \cite{B2}, we use the tools developed here to show that their centers are also 3-dimensional singularities, but with the strange property that they contain positive dimensional `smeared-out' points.  
Furthermore, in \cite{B4} we describe how this nonlocal geometry is reflected in the homology of their noncommutative residue fields.

(2) \textit{Quantum entanglement and nonlocality.}
In quantum physics, two entangled particles can effect each instantaneously no matter how great their spatial separation.
In \cite{B5,B6}, we introduce a toy model where such entangled states are supported on nonnoetherian singularities of spacetime.
Furthermore, we show that the noncommutative blowup of the singular support of an entangled Bell state leads to the notion that the \textit{nonlocal commutative} spacetime that we observe emerges from an underlying \textit{local noncommutative} spacetime \cite[Theorem C]{B5}.

(3) \textit{Nonnoetherian geometry.}  
In algebraic geometry, one may ask what is the geometry of a nonnoetherian algebra $R$?
The usual answer is the affine scheme $\operatorname{Spec}R$ whose global sections are $R$.
However, schemes are abstract objects, and we would like an answer to this question that we can visualize.
We aim to give such an answer here.

We now outline our main results.

In Section \ref{Section 1} we introduce the idea of a nonlocal algebraic variety.
To illustrate this notion, 
consider the algebra $S = k[x,y]$ and its nonnoetherian subalgebra
\begin{equation} \label{R =}
R = k\left[x,xy,xy^2, \ldots, \right] = k + xS.
\end{equation}
The maximal ideal spectrum $\operatorname{Max}R$ of $R$ may be viewed as 2-dimensional affine space $\mathbb{A}^2_k = \operatorname{Max}S$ with the line
$$\mathcal{Z}(x) = \left\{ x = 0 \right\} \subset \mathbb{A}^2_k$$
identified as a single closed point.
From this perspective, $\mathcal{Z}(x)$ is a 1-dimensional `smeared-out' point of $R$, and therefore $\operatorname{Max}R$ is nonlocal.

More generally, let $S$ be an integral domain and a finitely generated $k$-algebra, and let $R$ be a (possibly nonnoetherian) subalgebra of $S$.
In order to capture the locus where $\operatorname{Max}R$ `looks like' the variety $\operatorname{Max}S$, we introduce the open subset
$$U_{S/R} := \left\{ \mathfrak{n} \in \operatorname{Max}S \ | \ R_{\mathfrak{n}\cap R} = S_{\mathfrak{n}} \right\}.$$ 

\begin{Theorem A} (Theorem \ref{U and W}.)
Suppose $U_{S/R} \not = \emptyset$.
Then $\operatorname{Max}S$ and $\operatorname{Max}R$ are isomorphic on open dense subsets, and thus birationally equivalent. 
Furthermore, the Krull dimensions of $R$ and $S$ are equal.
\end{Theorem A}

In example (\ref{R =}), $U_{S/R}$ is the complement to the subvariety $\mathcal{Z}(x)$. 
We generalize this example by showing that if $R$ is generated by a subalgebra of $S$ and an ideal $I \subset S$, then $U_{S/R}$ contains the complement to the subvariety $\mathcal{Z}(I)$ in $\operatorname{Max}S$.
In addition, if $I$ is a non-maximal radical ideal of $S$ and $R = k + I$, then $U_{S/R} = \mathcal{Z}(I)^c$ (Proposition \ref{I}).

To formalize these notions, we introduce the following definitions.

\begin{Definition B} (Definition \ref{dimension def}.) \rm{ \ 
\begin{itemize}
 \item A finitely generated $k$-algebra $S$ is a \textit{depiction} of a subalgebra $R \subseteq S$ if 
\begin{enumerate}[(i)] 
 \item the morphism $\iota_S: \operatorname{Spec}S \rightarrow \operatorname{Spec}R$, $\mathfrak{q} \mapsto \mathfrak{q} \cap R$, is surjective,
 \item for each $\mathfrak{n} \in \operatorname{Max}S$, $R_{\mathfrak{n} \cap R}$ is noetherian iff $R_{\mathfrak{n} \cap R} = S_{\mathfrak{n}}$, and
 \item $U_{S/R} \not = \emptyset$.
\end{enumerate}  
 \item The \textit{geometric codimension} or \textit{geometric height} of $\mathfrak{p} \in \operatorname{Spec}R$ is the infimum 
$$\operatorname{ght}(\mathfrak{p}) := \operatorname{inf} \left\{ \operatorname{ht}(\mathfrak{q}) \ | \ \mathfrak{q} \in \iota^{-1}_S(\mathfrak{p}), \ S \text{ a depiction of } R \right\}.$$
The \textit{geometric dimension} of $\mathfrak{p}$ is the difference
$$\operatorname{gdim} \mathfrak{p} := \operatorname{dim}R - \operatorname{ght}(\mathfrak{p}).$$
\end{itemize}
} \end{Definition B}

In example (\ref{R =}), $R$ is depicted by $S$, and the geometric dimension of the closed point $\mathcal{Z}(x)$ of $R$ is 1 (Example \ref{ho ho ho}).
The following theorem describes the fundamental geometry of nonnoetherian algebras with finite Krull dimension.

\begin{Theorem C} (Theorems \ref{codim theorem}, \ref{non-max}, \ref{pos dim theorem}, and Proposition \ref{isolated sing}.)
Suppose $R$ admits a depiction $S$ and let $\mathfrak{p} \in \operatorname{Spec}R$.
Then
$$\operatorname{ght}(\mathfrak{p}) \leq \operatorname{ht}(\mathfrak{p}),$$
with equality if there is some $\mathfrak{q} \in \operatorname{Spec}S$ for which $\mathfrak{q} \cap R = \mathfrak{p}$ and $\mathcal{Z}\left( \mathfrak{q} \right) \cap U_{S/R} \not = \emptyset$.

Furthermore, the following are equivalent:
\begin{enumerate} 
  \item $R$ is noetherian.
  \item $U_{S/R} = \operatorname{Max}S$.
  \item $R = S$.
\end{enumerate}
In particular, if $R$ is noetherian, then its only depiction is itself.

Finally, if $\operatorname{Max}R$ contains a closed point of positive geometric dimension, then $R$ is nonnoetherian.
Conversely, if $R$ is an isolated nonnoetherian singularity, then $\operatorname{Max}R$ contains a closed point of positive geometric dimension.
\end{Theorem C}

Consequently, if $I$ is a radical ideal of a finitely generated $k$-algebra $S$, then the ring $R = k + I$ will be nonnoetherian if and only if $\operatorname{dim}\mathcal{Z}(I) \geq 1$ (Corollary \ref{R = k + I}).

We conclude the section by showing that depictions which are minimal with respect to inclusion do not exist in general, and maximal depictions are not unique in general (Proposition \ref{min and max}).

In Section \ref{depictions from impressions}, we study nonlocality in the context of noncommutative algebraic geometry.
Let $A$ be a finitely generated noncommutative $k$-algebra with center $Z$. 
We consider algebras with the following particularly nice matrix ring structure.

\begin{Definition} \label{impression definition} \rm{
An \textit{impression} of $A$ is a finitely generated commutative $k$-algebra $B$ and an algebra monomorphism $\tau: A \hookrightarrow M_d(B)$ such that (i) for each $\mathfrak{b}$ in some open dense subset of $\operatorname{Max}B$, the composition with the evaluation map
$$A \stackrel{\tau}{\longrightarrow} M_d(B) \stackrel{\epsilon_{\mathfrak{b}}}{\longrightarrow} M_d\left(B/\mathfrak{b} \right)$$
is surjective, and (ii) the morphism $\operatorname{Max}B \rightarrow \operatorname{Max}\tau(Z)$, $\mathfrak{b} \mapsto \mathfrak{b}1_d \cap \tau(Z)$, is surjective \cite[Definition and Lemma 2.1]{B}.
} \end{Definition}

An impression is useful in part because it determines the center $Z$ of $A$ \cite[Lemma 2.1]{B}.
Furthermore, if $A$ is a finitely generated $Z$-module, then $\tau$ determines all simple $A$-module isoclasses of maximal $k$-dimension \cite[Proposition 2.5]{B}.

Denote by $E_{ji} \in M_d(k)$ the matrix with a 1 in the $(ji)$-th slot and zeros elsewhere.
Let $A = kQ/I$ be a quiver algebra with vertex set $Q_0 = \left\{ 1, \ldots, d \right\}$, and suppose $A$ admits an impression $\tau: A \hookrightarrow M_d(B)$ such that $\tau(e_i) = E_{ii}$ for each $i \in Q_0$. 
For $p \in e_jAe_i$, denote by $\bar{\tau}(p) \in B$ the single non-zero entry of the matrix $\tau(p)$; whence
$$\tau(p) = \bar{\tau}(p)E_{ji}.$$
For each $i,j \in Q_0$, $\bar{\tau}$ defines a $k$-linear map $\bar{\tau}: e_jAe_i \rightarrow B$. 
Set
$$\begin{array}{c}
R := k\left[ \displaystyle{\cap_{i \in Q_0}} \bar{\tau}\left( e_iAe_i \right) \right] \subseteq B,\\
S := k \left[ \cup_{i \in Q_0} \bar{\tau}\left( e_iAe_i \right) \right] \subseteq B.
\end{array}$$

\begin{Theorem D} (Theorem \ref{impression corollary}.)
Suppose $\tau: A \hookrightarrow M_{|Q_0|}(B)$ is an impression of $A$ with $B$ an integral domain and $\tau(e_i) = E_{ii}$ for each $i \in Q_0$.
Then $U_{S/R} \not = \emptyset$.
Furthermore, if $\left\{ \mathfrak{n} \in \operatorname{Max}S \ | \ R_{\mathfrak{n} \cap R} \text{ is noetherian} \right\} \subseteq U_{S/R}$, then
\begin{enumerate}
 \item The center $Z$ of $A$ is isomorphic to $R$ and is depicted by $S$.
 \item The statements
  \begin{enumerate} 
   \item $R = S$.
   \item $A$ is a finitely generated $Z$-module.
   \item $Z$ is noetherian.
   \item $A$ is noetherian.
  \end{enumerate}
  satisfy the equivalences
  $$\xy 
 (-24,0)*+{\text{(a)}}="1";(-8,0)*+{\text{(b)}}="2";(8,0)*+{\text{(c)}}="3";(24,0)*+{\text{(d)}}="4";
{\ar@{<=>}"1";"2"};{\ar@{<=>}"2";"3"};{\ar@{=>}@/^/"3";"4"};{\ar@{=>}@/^/^{\star}"4";"3"};
  \endxy$$
  where $(\star)$ holds if the $\bar{\tau}$-image of a path is a monomial in $B$ and $I$ is generated by binomials in the paths of $Q$.
\end{enumerate}
\end{Theorem D}

Again consider the nonnoetherian algebra $R = k + xS$ in example (\ref{R =}).
In Example \ref{warp drive} and Proposition \ref{A quiver}, we study the endomorphism ring 
$$A = \operatorname{End}_R\left(R \oplus xS \right)$$
of the reflexive $R$-module $R \oplus xS$.  
This algebra may be viewed as the noncommutative blowup of $R$ at the isolated singular point $xS$ of $\operatorname{Max}R$ \cite[Section R]{L}. 
Furthermore, $A$ admits an impression $\tau: A \hookrightarrow M_2(S)$ satisfying 
$$R = k + xS = k\left[ \displaystyle{\cap_{i \in Q_0}} \bar{\tau}\left( e_iAe_i \right) \right] \ \ \text{ and } \ \ S = k[x,y] = k\left[ \displaystyle{\cup_{i \in Q_0}} \bar{\tau}\left( e_iAe_i \right) \right].$$
By Theorem D, $A$ is nonnoetherian, has center $R$, and is an infinitely generated $R$-module.
We show that the projective dimension and geometric codimension of its vertex simple modules $V_i$ coincide:
$$\operatorname{pd}_A\left(V_i \right) = \operatorname{ght}\left( \operatorname{ann}_{e_iAe_i}\left( V_i \right) \right).$$
This example suggests that the notion of geometric dimension is, in a suitable sense, a natural definition. 

\textit{Notation:}
We will denote by $\operatorname{dim}R$ the Krull dimension of $R$; by $\operatorname{Frac}R$ the ring of fractions of $R$; by $\operatorname{Max}R$ the set of maximal ideals of $R$; and by $\operatorname{Spec}R$ either the set of prime ideals of $R$ or the affine $k$-scheme with global sections $R$.  
For $\mathfrak{a} \subset R$ we will denote by $\mathcal{Z}(\mathfrak{a})$ either the closed set $\left\{ \mathfrak{m} \in \operatorname{Max}R \ | \ \mathfrak{m} \supseteq \mathfrak{a} \right\}$ of $\operatorname{Max}R$ or the closed set $\left\{ \mathfrak{p} \in \operatorname{Spec}R \ | \ \mathfrak{p} \supseteq \mathfrak{a} \right\}$ of $\operatorname{Spec}R$, depending on the context.
For a subset $Y$ of $\operatorname{Max}S$, set $Y^c := \operatorname{Max}S \setminus Y$.

In Section \ref{depictions from impressions} we will denote by $Q = \left( Q_0,Q_1,\operatorname{t}, \operatorname{h} \right)$ a quiver with vertex set $Q_0$, arrow set $Q_1$, and head and tail maps $\operatorname{h},\operatorname{t}: Q_1 \to Q_0$.
We will denote by $kQ$ the path algebra of $Q$, and by $e_i$ the idempotent corresponding to the vertex $i \in Q_0$.
Multiplication of paths is read right to left, following the composition of maps.
By module we mean left module.
By infinitely generated $R$-module, we mean a module that is not finitely generated over $R$.

\section{Nonnoetherian geometry as nonlocal geometry} \label{Section 1}

Throughout $k$ is an algebraically closed field; $S$ is an integral domain and a noetherian $k$-algebra; and $R$ is a (possibly nonnoetherian) subalgebra of $S$. 
We begin with the following well known lemma.

\begin{Lemma} \label{subset}
If $\mathfrak{q} \in \operatorname{Spec}S$, then $\mathfrak{q} \cap R \in \operatorname{Spec}R$.
Furthermore, if $S$ is a finitely generated $k$-algebra and $\mathfrak{n} \in \operatorname{Max}S$, then $\mathfrak{n} \cap R \in \operatorname{Max}R$.
\end{Lemma}

\begin{proof}
We show the second statement. 
Suppose $S$ is a finitely generated $k$-algebra and let $\mathfrak{n} \in \operatorname{Max}S$.
Then $S/\mathfrak{n} \cong k$ since $S$ is finitely generated over the algebraically closed field $k$.\footnote{This statement is false in general without the assumption that $S$ is a finitely generated $k$-algebra.
Indeed, let $S$ be $\mathbb{C}$, let $k$ be the algebraic closure of $\mathbb{Q}$, and let $\mathfrak{n}$ be the maximal ideal $0$ of $\mathbb{C}$.
Then $S/\mathfrak{n} \not \cong k$.}
Thus the composition $\psi: R \hookrightarrow S \rightarrow S/\mathfrak{n}$ is surjective since $1_S \in R$.  
Whence $R/ \operatorname{ker}\psi \cong k$. 
Therefore $\operatorname{ker}\psi = \mathfrak{n} \cap R$ is a maximal ideal of $R$.
\end{proof}

Consider the morphisms
$$\begin{array}{lcccclccc}
\iota: & \operatorname{Spec}S & \to & \operatorname{Spec}R & \ \ \ \text{ and } \ \ \ &  \kappa: & \operatorname{Max}S & \to & \operatorname{Spec}R \\
& \mathfrak{q} & \mapsto & \mathfrak{q} \cap R & & & \mathfrak{n} &\mapsto & \mathfrak{n} \cap R.
\end{array}$$
Let $\mathfrak{q} \in \operatorname{Spec}S$ and set $\mathfrak{p} = \iota(\mathfrak{q}) = \mathfrak{q} \cap R$.
Then $R \setminus \mathfrak{p} \subseteq S \setminus \mathfrak{q}$.
Whence
\begin{equation} \label{Rp in Sq}
R_{\iota(\mathfrak{q})} = R_{\mathfrak{p}} \subseteq S_{\mathfrak{q}}.
\end{equation}
Thus the embedding $R \hookrightarrow S$ induces the morphism of schemes 
$$\iota: \left(\operatorname{Spec}S, \mathcal{O}_{\operatorname{Spec}S}\right) \longrightarrow \left( \operatorname{Spec}R, \mathcal{O}_{\operatorname{Spec}R} \right).$$

To aid our analysis of $R$, we introduce the following subsets of $\operatorname{Max}S$ and $\operatorname{Spec}S$:
\begin{equation} \label{U sets}
\begin{array}{l}
U_{S/R} := \left\{ \mathfrak{n} \in \operatorname{Max}S \ | \ R_{\mathfrak{n} \cap R} = S_{\mathfrak{n}} \right\},\\
\tilde{U}_{S/R} := \left\{ \mathfrak{q} \in \operatorname{Spec}S \ | \ R_{\mathfrak{q} \cap R} = S_{\mathfrak{q}} \right\},\\
U^*_{S/R} := \left\{ \mathfrak{n} \in \operatorname{Max}S \ | \ R_{\mathfrak{n} \cap R} \text{ is noetherian} \right\},\\
W_{S/R} := \left\{ \mathfrak{n} \in \operatorname{Max}S \ | \ \sqrt{(\mathfrak{n}\cap R)S} = \mathfrak{n} \right\}.
\end{array}
\end{equation}
These subsets will play a central role throughout this paper.
Furthermore, $U_{S/R}$ will play a central role in its sequels \cite{B2,B3} in the context of dimer algebras. 
If $R$ and $S$ are fixed, then we will often omit the subscript $S/R$.

\begin{Lemma} \label{Rp Sq}
Suppose $\mathfrak{p} \in \operatorname{Spec}R$ and $\mathfrak{q} \in \iota^{-1}(\mathfrak{p})$. 
If $\mathcal{Z}(\mathfrak{q}) \cap U \not = \emptyset$, then $R_{\mathfrak{p}} = S_{\mathfrak{q}}$.
\end{Lemma}

\begin{proof}
Let $\mathfrak{n} \in \mathcal{Z}(\mathfrak{q}) \cap U$ and set $\mathfrak{m} = \mathfrak{n} \cap R$.
Then $R_{\mathfrak{m}} = S_{\mathfrak{n}}$.
Furthermore, $R_{\mathfrak{p}} \subseteq S_{\mathfrak{q}}$ by (\ref{Rp in Sq}).
Thus, since $\mathfrak{p} \subseteq \mathfrak{q}$,
$$S_{\mathfrak{q}} = \left( S_{\mathfrak{n}} \right)_{\mathfrak{q}} \subseteq \left( S_{\mathfrak{n}} \right)_{\mathfrak{p}} = \left( R_{\mathfrak{m}} \right)_{\mathfrak{p}} = R_{\mathfrak{p}} \subseteq S_{\mathfrak{q}}.$$
Therefore $R_{\mathfrak{p}} = S_{\mathfrak{q}}$.
\end{proof}

\begin{Lemma} \label{qq'}
Suppose $\mathfrak{q}, \mathfrak{q}' \in \operatorname{Spec}S$ satisfy 
$$\mathfrak{q} \cap R = \mathfrak{q}' \cap R \ \ \text{ and } \ \ \mathfrak{q} \subseteq \mathfrak{q}'.$$
If $\mathcal{Z}(\mathfrak{q}) \cap U \not = \emptyset$ or $\mathcal{Z}(\mathfrak{q}') \cap U \not = \emptyset$, then $\mathfrak{q} = \mathfrak{q}'$.
\end{Lemma}

\begin{proof}
We claim that $S_{\mathfrak{q}} = S_{\mathfrak{q}'}$, and consequently $\mathfrak{q} = \mathfrak{q}'$ since $S_{\mathfrak{q}}$ has a unique maximal ideal.

Indeed, suppose $\mathcal{Z}(\mathfrak{q}) \cap U \not = \emptyset$.
Then by Lemma \ref{Rp Sq},
$$S_{\mathfrak{q}'} \subseteq S_{\mathfrak{q}} = R_{\mathfrak{q}\cap R} = R_{\mathfrak{q}' \cap R} \subseteq S_{\mathfrak{q}'}.$$
Therefore $S_{\mathfrak{q}} = S_{\mathfrak{q}'}$.

So suppose $\mathcal{Z}(\mathfrak{q}') \cap U \not = \emptyset$. 
Since $\mathfrak{q} \subseteq \mathfrak{q}'$, we have 
$$\mathcal{Z}(\mathfrak{q}') \subseteq \mathcal{Z}(\mathfrak{q}).$$
Whence $\mathcal{Z}(\mathfrak{q}) \cap U \not = \emptyset$, which was the previous case.
\end{proof}

In the remainder of this section we assume that $S$ is a finitely generated $k$-algebra, unless stated otherwise.
Recall that $S$ is an overring of a domain $R$ if $R \subseteq S \subseteq \operatorname{Frac}R$.
\newpage
\begin{Proposition} \label{electric}
Suppose $U$ is nonempty. 
Then 
\begin{enumerate}
 \item $S$ is an overring of $R$; and
 \item $U$ is an open dense subset of $\operatorname{Max}S$.
\end{enumerate}
In particular, the function fields of $\operatorname{Spec}S$ and $\operatorname{Spec}R$ coincide. 
\end{Proposition}

\begin{proof}
(1) Suppose $\mathfrak{n} \in U := U_{S/R}$.
Since $S$ is an integral domain, we have
\begin{equation} \label{frac S}
\operatorname{Frac}S = \operatorname{Frac}(S_{\mathfrak{n}}) = \operatorname{Frac}(R_{\mathfrak{n} \cap R}) =\operatorname{Frac}R.
\end{equation}

(2) We first claim that $U$ coincides with the locus
$$U' := \left\{ \mathfrak{n} \in \operatorname{Max}S \ | \ S \subseteq R_{\mathfrak{n} \cap R} \right\}.$$
Indeed, let $\mathfrak{n} \in U$.
Then $S_{\mathfrak{n}} = R_{\mathfrak{n} \cap R}$.
Whence $S \subseteq R_{\mathfrak{n} \cap R}$.
Thus $U \subseteq U'$.

Conversely, suppose $\mathfrak{n} \in U'$.
Then 
$$S_{\mathfrak{n}} = S(S\setminus \mathfrak{n} )^{-1} \subseteq R_{\mathfrak{n} \cap R}(R_{\mathfrak{n} \cap R} \setminus (\mathfrak{n} \cap R_{\mathfrak{n} \cap R}))^{-1} = R_{\mathfrak{n} \cap R} = S_{\mathfrak{n}}.$$
Whence $R_{\mathfrak{n} \cap R} = S_{\mathfrak{n}}$.
Thus $U \supseteq U'$.
Therefore $U = U'$, proving our claim.

Suppose that $U = U'$ is nonempty.
By assumption, $S$ is a finitely generated $k$-algebra.
Thus there is a finite set $T := \left\{ s_1, \ldots, s_{\ell} \right\} \subset S$, minimal with respect to inclusion, such that
$$S = R\left[ s_1, \ldots, s_{\ell} \right].$$
By Claim (1), there are elements $r_1, \ldots, r_{\ell} \in R$ such that $r_i s_i \in R$.
Whence $s_i \in R\left[ r_i^{-1} \right]$.
Therefore
$$S \subseteq R\left[ r_1^{-1}, \ldots, r_{\ell}^{-1} \right].$$

Consider the open set
$$U''_T := \bigcap_{1 \leq i \leq \ell} D(r_i),$$
where $D(r_i) = \left\{ \mathfrak{n} \in \operatorname{Max}S \ | \ \mathfrak{n} \not \ni r_i \right\}$ is the principal open set where $r_i$ does not vanish.
Then
$$U''_T \subseteq U'.$$
Furthermore, $U'$ is the union over all such minimal sets $T \subset S$,
$$U' = \bigcup_{T} U''_T.$$
In particular, $U'$ is open.
Thus $U'$ is dense since $S$ is an integral domain and $U'$ is nonempty.
Therefore $U = U'$ is open dense.
\end{proof} 

The morphism $\iota: \operatorname{Spec}S \to \operatorname{Spec}R$ is injective on the subset $\tilde{U} \subseteq \operatorname{Spec}S$ defined in (\ref{U sets}).
Indeed, if $\mathfrak{q} \cap R = \mathfrak{q}' \cap R$ and $\mathfrak{q}, \mathfrak{q}' \in \tilde{U}$, then
$$S_{\mathfrak{q}} = R_{\mathfrak{q} \cap R} = R_{\mathfrak{q}' \cap R} = S_{\mathfrak{q}'}.$$
Whence $\mathfrak{q} = \mathfrak{q}'$.
This fact is generalized in the following theorem.

Furthermore, we will find that the locus $W \subset \operatorname{Max}S$ is similar in spirit to the Azumaya locus of $A$ when $A$ is a noncommutative algebra, module-finite over its center $Z$.\footnote{Recall that if $\mathfrak{n}, \mathfrak{n}' \in \operatorname{Max}A$ and $\mathfrak{n} \cap Z = \mathfrak{n}' \cap Z$ is in the Azumaya locus of $A$, then $\mathfrak{n} = \mathfrak{n}'$.}  

\begin{Theorem} \label{U and W} 
Let $R$ be a subalgebra of $S$.
\begin{enumerate}
 \item If $R_{\mathfrak{q}\cap R} = S_{\mathfrak{q}}$, then 
$$\iota^{-1}\iota(\mathfrak{q}) = \left\{ \mathfrak{q} \right\}.$$
In particular, this holds if $\mathcal{Z}(\mathfrak{q}) \cap U_{S/R} \not = \emptyset$.
 \item The locus $W_{S/R}$ is the subset of all $\mathfrak{n} \in \operatorname{Max}S$ for which
$$\kappa^{-1}\kappa(\mathfrak{n}) = \left\{ \mathfrak{n} \right\}.$$ 
In particular, $U_{S/R} \subseteq W_{S/R}$.  
 \item If $U_{S/R} \not = \emptyset$, then $\operatorname{Max}S$ and $\operatorname{Max}R$ are isomorphic on open dense subsets, and thus birationally equivalent.
 \item If $U_{S/R} \not = \emptyset$, then the Krull dimensions of $R$ and $S$ coincide,
 $$\operatorname{dim}R = \operatorname{dim}S.$$
\end{enumerate}
\end{Theorem}

\begin{proof}
(1.i) Let $\mathfrak{q} \in \tilde{U}_{S/R}$, and suppose $\mathfrak{q}' \in \iota^{-1}\iota(\mathfrak{q})$.  
We want to show that $\mathfrak{q}' = \mathfrak{q}$. 

Set $\mathfrak{p} := \mathfrak{q} \cap R = \mathfrak{q}' \cap R$.
Then
$$\mathfrak{p}_{\mathfrak{p}} = \left( \mathfrak{q} \cap R \right) R_{\mathfrak{p}} \subseteq \mathfrak{q}_{\mathfrak{q}} \cap R_{\mathfrak{p}}.$$
Thus, since $\mathfrak{p}_{\mathfrak{p}}$ is maximal in $R_{\mathfrak{p}}$, we have $\mathfrak{p}_{\mathfrak{p}} = \mathfrak{q}_{\mathfrak{q}} \cap R_{\mathfrak{p}}$.
Similarly $\mathfrak{p}_{\mathfrak{p}} = \mathfrak{q}'_{\mathfrak{q}} \cap R_{\mathfrak{p}}$.
Whence
\begin{equation} \label{1.i.3}
\mathfrak{q}_{\mathfrak{q}} \cap R_{\mathfrak{p}} = \mathfrak{p}_{\mathfrak{p}} = \mathfrak{q}'_{\mathfrak{q}} \cap R_{\mathfrak{p}}.
\end{equation}
In particular, since $\mathfrak{q} \in \tilde{U}_{S/R}$, 
$$U_{S_{\mathfrak{q}}/R_{\mathfrak{p}}} := \left\{ \mathfrak{n} \in \operatorname{Max}S_{\mathfrak{q}} \ | \ \left( R_{\mathfrak{p}} \right)_{\mathfrak{n} \cap R_{\mathfrak{p}}} = \left( S_{\mathfrak{q}} \right)_{\mathfrak{n}} \right\} = \left\{ \mathfrak{q}_{\mathfrak{q}} \right\}.$$

Furthermore, since $\mathfrak{q} \in \tilde{U}_{S/R}$ and $\mathfrak{q}' \in \iota^{-1}\iota(\mathfrak{p})$,
$$S_{\mathfrak{q}} = R_{\mathfrak{p}} \subseteq S_{\mathfrak{q}'}.$$ 
Thus $\mathfrak{q}' \subseteq \mathfrak{q}$.
Therefore
\begin{equation} \label{1.i.1}
\mathfrak{q}'_{\mathfrak{q}} \subseteq \mathfrak{q}_{\mathfrak{q}}.
\end{equation}
Whence
\begin{equation} \label{1.i.2}
\mathcal{Z}\left( \mathfrak{q}'_{\mathfrak{q}} \right) \cap U_{S_{\mathfrak{q}}/R_{\mathfrak{p}}} = \left\{ \mathfrak{q}_{\mathfrak{q}} \right\} \not = \emptyset.
\end{equation}

Thus, by Lemma \ref{qq'}, (\ref{1.i.3}), (\ref{1.i.1}), and (\ref{1.i.2}) imply that
$$\mathfrak{q}'_{\mathfrak{q}} = \mathfrak{q}_{\mathfrak{q}}.$$
Therefore $\mathfrak{q}' = \mathfrak{q}'_{\mathfrak{q}} \cap S = \mathfrak{q}_{\mathfrak{q}} \cap S = \mathfrak{q}$, which is what we wanted to show.

(1.ii) If $\mathcal{Z}(\mathfrak{q}) \cap U_{S/R} \not = \emptyset$, then $\mathfrak{q} \in \tilde{U}_{S/R}$ by Lemma \ref{Rp Sq}.

(2) We now claim that $\mathfrak{n} \in W^c_{S/R}$ if and only if there is a point $\mathfrak{n}' \in \operatorname{Max}S$, not equal to $\mathfrak{n}$, such that $\kappa(\mathfrak{n}) = \kappa(\mathfrak{n}')$.  

First note that for any $\mathfrak{m} \in \operatorname{Max}R$ and $\mathfrak{n} \in \kappa^{-1}(\mathfrak{m})$, we have $\mathfrak{m} \subseteq \mathfrak{m}S \subseteq \sqrt{\mathfrak{m}S} \subseteq \mathfrak{n}$.
Thus $\mathfrak{m} \subseteq \sqrt{\mathfrak{m}S} \cap R \subseteq \mathfrak{n} \cap R = \mathfrak{m}$. 
Whence
\begin{equation} \label{expanding}
\sqrt{\mathfrak{m}S} \cap R = \mathfrak{m}.
\end{equation}

Now let $\mathfrak{n} \in W^c_{S/R}$ and set $\mathfrak{m}:= \mathfrak{n} \cap R$. 
Then by definition, $\mathfrak{n} \not = \sqrt{\mathfrak{m}S}$.  
Since $S$ is Jacobson we have 
$$\sqrt{\mathfrak{m}S} = \ \ \ \ \bigcap_{\mathclap{\mathfrak{m}S \subseteq \mathfrak{q} \in \operatorname{Max}S}}\mathfrak{q}.$$
Thus there exists a maximal ideal $\mathfrak{n}' \not = \mathfrak{n}$ of $S$ such that $\sqrt{\mathfrak{m}S} \subseteq \mathfrak{n}'$.  
Whence
$$\kappa(\mathfrak{n}) = \mathfrak{n} \cap R = \mathfrak{m} \stackrel{(\ref{expanding})}{=} \sqrt{\mathfrak{m}S} \cap R \subseteq \mathfrak{n}'\cap R = \kappa(\mathfrak{n}').$$
Therefore $\kappa(\mathfrak{n}) = \kappa(\mathfrak{n}')$ by Lemma \ref{subset}.

Conversely, suppose there are distinct points $\mathfrak{n}, \mathfrak{n}' \in \operatorname{Max}S$ such that $\kappa(\mathfrak{n}) = \kappa(\mathfrak{n}')$.  
Then $\mathfrak{n} \cap R = \mathfrak{n}' \cap R =: \mathfrak{m}$.
Therefore 
$$\sqrt{\mathfrak{m}S} \subseteq \mathfrak{n} \cap \mathfrak{n}' \subsetneq \mathfrak{n}.$$

(3) Suppose $U_{S/R} \not = \emptyset$.
By Claim (1) $\iota$ is injective on $U_{S/R}$.
Furthermore, $U_{S/R}$ is an open dense subset of $\operatorname{Max}S$ by Proposition \ref{electric}.
Therefore $\operatorname{Max}S$ and $\operatorname{Max}R$ are birationally equivalent.

(4) Finally, we show that $\operatorname{dim}R = \operatorname{dim}S$.
Fix $\mathfrak{n} \in U_{S/R}$ and set $\mathfrak{m} := \mathfrak{n} \cap R$.
Then
$$\operatorname{dim}R \leq \operatorname{trdeg}_k \operatorname{Frac}R  \stackrel{\textsc{(i)}}{=} \operatorname{trdeg}_k \operatorname{Frac}S \stackrel{\textsc{(ii)}}{=} \operatorname{dim}S \stackrel{\textsc{(iii)}}{=} \operatorname{dim}S_{\mathfrak{n}} \stackrel{\textsc{(iv)}}{=} \operatorname{dim}R_{\mathfrak{m}} \leq \operatorname{dim}R,$$
where (\textsc{i}) holds by Proposition \ref{electric}; (\textsc{ii}) and (\textsc{iii}) hold since $S$ is a noetherian integral domain over $k$; and (\textsc{iv}) holds since $\mathfrak{n}$ is in $U_{S/R}$.
\end{proof}

\begin{Example} \label{first ex} \rm{
Let $S = k[x,y]$ and $R = k[x,xy,xy^2, \ldots] = k + xS$.  
For any $b \in k$, the ideals $(x,y-b)S, xS \in \operatorname{Spec}S$ satisfy
$$(x,y-b)S \cap R = xS \cap R = (x,xy,xy^2,\ldots) \in \operatorname{Max}R.$$
Thus $(x,y-b)S \in W^c$ by Theorem \ref{U and W}.2.
}\end{Example}

\begin{Remark} \label{U not = W} \rm{
In general, $U$ need not equal $W$.
Indeed, consider the algebras
$$S = k[x] \ \ \text{ and } \ \ R = k + x^2S = k[x^2,x^3] \cong k[u,v]/(u^3 - v^2).$$
Then $U = \mathbb{A}^1 \setminus \{ 0 \}$ and $W = \mathbb{A}^1$.
}\end{Remark}

The following proposition generalizes the fact that if $R$ is a finitely generated $k$-algebra and $\mathfrak{m} \in \operatorname{Max}R$, then
$$R = k + \mathfrak{m}.$$
Conversely, if $I$ is an ideal in $S$ and $R = k + I$, then $\mathcal{Z}(I)$ is a closed point in $\operatorname{Spec}R$.
In Corollary \ref{R = k + I} below we will show that $R = k + I$ is nonnoetherian whenever the dimension of the subvariety $\mathcal{Z}(I)$ in $\operatorname{Max}S$ is nonzero.

\begin{Proposition} \label{I}
Let $S$ be an integral domain and a noetherian $k$-algebra.
Consider a subalgebra $R'$ of $S$, an ideal $I$ of $S$, and form the algebra 
$$R = k[R', I].$$
Then $U$ (hence $W$) contains the open subset $\mathcal{Z}(I)^c$ of $\operatorname{Max}S$.  

Furthermore, if $I \subset S$ is a non-maximal radical ideal and
$$R = k[I] = k+I,$$
then $U = W = \mathcal{Z}(I)^c$.
\end{Proposition}

\begin{proof}
First suppose $R = k[R',I]$.
We claim that if $\mathfrak{q} \in \operatorname{Spec}S$ does not contain $I$, then $R_{\mathfrak{q} \cap R} = S_{\mathfrak{q}}$; in particular, 
\begin{equation} \label{ZIc}
\mathcal{Z}(I)^c \subseteq U.
\end{equation}

Set $\mathfrak{p}:= \mathfrak{q} \cap R$.  
Then $R_{\mathfrak{p}} \subseteq S_{\mathfrak{q}}$.
To show the reverse inclusion, suppose $a \in S_{\mathfrak{q}}$.
Then there is some $f,g \in S$ with $g \not \in \mathfrak{q}$ such that $a = \frac fg$.  
Since $\mathfrak{q}$ does not contain $I$, there is some $c \in I \setminus \mathfrak{q}$.  
Furthermore, since $c,g \in S \setminus \mathfrak{q}$ and $\mathfrak{q}$ is prime, we have 
$$cg \in S \setminus \mathfrak{q}.$$  
Since $c \in I$, we have $cg \in I \subset R$.
Whence 
$$cg \in R \setminus \mathfrak{p}.$$  
But also 
$$b := agc = fc \in I \subset R.$$ 
Thus $a = \frac{b}{cg} \in R_{\mathfrak{p}}$.
Therefore $S_{\mathfrak{q}} \subseteq R_{\mathfrak{p}}$.
Consequently $R_{\mathfrak{q} \cap R} = S_{\mathfrak{q}}$.

Now suppose $R = k[I]$, where $I$ is a non-maximal radical ideal of $S$.
Let $\mathfrak{n} \in \mathcal{Z}(I)$.  
Then $\mathfrak{n} \supseteq I$, so 
$$\mathfrak{n} \cap R \supseteq I \cap R = I.$$
Whence $\mathfrak{n} \cap R = I$ since $I$ is a maximal ideal of $R$.  
But $\sqrt{IS} = IS = I \subsetneq \mathfrak{n}$ since $I$ is a radical ideal of $S$.
Thus $\mathcal{Z}(I) \subseteq W^c$.  

Therefore
$$\mathcal{Z}(I) \subseteq W^c 
\stackrel{\textsc{(i)}}{\subseteq} U^c 
\stackrel{\textsc{(ii)}}{\subseteq} \mathcal{Z}(I),$$
where (\textsc{i}) holds by Theorem \ref{U and W}.1 and (\textsc{ii}) holds by (\ref{ZIc}).
\end{proof}

\begin{Remark} \rm{
$U$ may properly contain $\mathcal{Z}(I)^c$; for example, take $R' = S$, in which case $U = \operatorname{Max}S$.  
}\end{Remark}

\begin{Example} \label{second ex} \rm{
\textit{A geometric picture.}

(i) Let $S = k[x,y]$ and $R = k + xS$.  
By Proposition \ref{I}, we can form $\operatorname{Max}R$ from $\operatorname{Max}S = \mathbb{A}^2_k$ by declaring the line 
$$\mathcal{Z}(x) = \{ x=0 \} \subset \mathbb{A}^2_k$$ 
to be a single (closed) point, and all other points, $U = \{ x \not = 0\}$, remain unaltered.
In this description, $\mathcal{Z}(x)$ appears to be a 1-dimensional, hence nonlocal, point of $\operatorname{Max}R$.

(ii) Let $S = k[x,y,z]$ and $R = k[x,y,yz,yz^2, \ldots] = k[x, yS]$.  
We can form $\operatorname{Max}R$ from $\operatorname{Max}S = \mathbb{A}^3_k$ by declaring each line 
$$\mathcal{Z}(x-c,y) = \{x=c, y=0\} \subset \mathbb{A}^3_k$$ 
to be a single point, and all other points, $U = \{y \not = 0\}$, remain unaltered.
Similar to the previous example, each subvariety $\mathcal{Z}(x-c,y)$ appears to be a 1-dimensional, hence nonlocal, point of $\operatorname{Max}R$.
} \end{Example}

\section{Depictions and the geometric dimension of a point}

Throughout $S$ is an integral domain and a finitely generated $k$-algebra.

We introduce the following definitions with the aim of constructing a geometric theory of nonnoetherian algebras, and in particular to formalize the geometric pictures in Example \ref{second ex}.
Recall that if $R$ is an integral domain and finitely generated over $k$, then the dimension of a point $\mathfrak{p} \in \operatorname{Spec}R$ coincides with the Krull dimension of $R$ minus the height of $\mathfrak{p}$,
$$\operatorname{dim} \mathfrak{p} := \operatorname{dim}R/\mathfrak{p} = \operatorname{dim}R - \operatorname{ht}(\mathfrak{p}).$$
The dimension of $\mathfrak{p}$ is then zero whenever $\mathfrak{p}$ is a maximal ideal.

\begin{Definition} \label{dimension def} \rm{ \ 
\begin{itemize}
 \item A finitely generated $k$-algebra $S$ is a \textit{depiction} of a subalgebra $R \subseteq S$ if
\begin{enumerate}[(i)] 
 \item the morphism $\iota_S: \operatorname{Spec}S \rightarrow \operatorname{Spec}R$, $\mathfrak{q} \mapsto \mathfrak{q} \cap R$, is surjective,
 \item for each $\mathfrak{n} \in \operatorname{Max}S$, $R_{\mathfrak{n} \cap R}$ is noetherian iff $R_{\mathfrak{n} \cap R} = S_{\mathfrak{n}}$, and
 \item $U_{S/R} \not = \emptyset$.
\end{enumerate}  
 \item The \textit{geometric codimension} or \textit{geometric height} of $\mathfrak{p} \in \operatorname{Spec}R$ is the infimum 
$$\operatorname{ght}(\mathfrak{p}) := \operatorname{inf} \left\{ \operatorname{ht}(\mathfrak{q}) \ | \ \mathfrak{q} \in \iota^{-1}_S(\mathfrak{p}), \ S \text{ a depiction of } R \right\}.$$
The \textit{geometric dimension} of $\mathfrak{p}$ is the difference
$$\operatorname{gdim} \mathfrak{p} := \operatorname{dim}R - \operatorname{ght}(\mathfrak{p}).$$
\end{itemize}
} \end{Definition}

\begin{Remark} \rm{
Note that condition (ii) is equivalent to $U^*_{S/R} \subseteq U_{S/R}$.
Furthermore, this condition is independent of conditions (i) and (iii).
Indeed, consider 
$$S = k[x,y], \ \ \ R = k[x,xy,y^2,y^3], \ \ \text{ and } \ \ \mathfrak{n} = (x,y)S \in \operatorname{Max}S.$$
Then conditions (i) and (iii) hold.\footnote{Naively it appears as though condition (i) may not hold since $xy \in xS \cap R \setminus xR$, and so there is no $\mathfrak{q} \in \operatorname{Spec}S$ for which $\mathfrak{q} \cap R = xR$.
However, $xR$ is not prime in $R$ since $x \cdot y^3 = (xy)y^2$.}
However, $R_{\mathfrak{n} \cap R}$ is noetherian whereas $R_{\mathfrak{n} \cap R} \not = S_{\mathfrak{n}}$.
}\end{Remark}

\begin{Remark} \rm{
Let $S$ and $S'$ be depictions of $R$, and let $\mathfrak{p} \in \operatorname{Spec}R$.  
Then in general the infimums of heights of ideals in $\iota_S^{-1}\left(\mathfrak{p} \right)$ and $\iota_{S'}^{-1}\left( \mathfrak{p} \right)$ do not coincide. For example, let 
$$S = k[x,y,z], \ \ \ S' = S[x^{-1}], \ \ \ R = k + x(y,z)S, \ \ \text{ and } \ \ \mathfrak{m} = x(y,z)S \in \operatorname{Max}R.$$
Then
$$\operatorname{inf} \left\{ \operatorname{ht}(\mathfrak{q}) \ | \ \mathfrak{q} \in \iota^{-1}_S\left( \mathfrak{m} \right) \right\} = \operatorname{ht}(xS) = 1,$$
whereas
$$\operatorname{inf} \left\{ \operatorname{ht}(\mathfrak{q}) \ | \ \mathfrak{q} \in \iota^{-1}_{S'}\left( \mathfrak{m} \right) \right\} = \operatorname{ht}((y,z)S') = 2.$$
}\end{Remark}

\begin{Remark} \rm{
In \cite{S}, Schwede gives a geometric description of subalgebras of noetherian algebras which is based on the gluing of schemes.  
In particular, the subalgebra $k[x,xy,xy^2,\ldots] \subset k[x,y]$ is described as the fiber product $k[x,y] \times_{k[y]} k$ (\cite[Example 3.7]{S}).
}\end{Remark}

\begin{Lemma} \label{nonlocal points}
If $S$ and $S'$ are depictions of $R$, then
$$\iota_S( U_{S/R}) = \iota_{S'}( U_{S'/R}).$$
\end{Lemma}

\begin{proof}
We claim that 
\begin{equation} \label{iota claim}
\iota_S(U^*_{S/R}) = U^*_{R/R}.
\end{equation}
Indeed, suppose $\mathfrak{m} \in U^*_{R/R}$.
Then $R_{\mathfrak{m}}$ is noetherian. 
Thus, since $\iota_S$ is surjective, $\mathfrak{m} \in \iota_S(U^*_{S/R})$. 
Conversely, suppose $\mathfrak{m} \in \iota_S(U^*_{S/R})$.
Then there is some $\mathfrak{n} \in \operatorname{Max}S$ for which $\mathfrak{n} \cap R = \mathfrak{m}$ and $R_{\mathfrak{n}\cap R}$ is noetherian.
Thus $\mathfrak{m} \in U^*_{R/R}$.
Therefore (\ref{iota claim}) holds.

Since $S$ and $S'$ are depictions of $R$, we have $U_{S/R} = U^*_{S/R}$ and $U_{S'/R} = U^*_{S'/R}$.
Therefore by (\ref{iota claim}),
$$\iota_S(U_{S/R}) = \iota_S(U^*_{S/R}) = U^*_{R/R} = \iota_{S'}(U^*_{S'/R}) = \iota_{S'}(U_{S'/R}).$$
\end{proof}

The following lemma will be useful in Section \ref{depictions from impressions}.

\begin{Lemma} \label{surjective1}
Let $R$ be a subalgebra of $S$, and suppose $k$ is uncountable.
Then the morphism $\iota: \operatorname{Spec}S \to \operatorname{Spec}R$ is surjective if and only if the morphism 
$\kappa: \operatorname{Max}S \to \operatorname{Max}R$ is surjective. 
\end{Lemma}

\begin{proof}
Suppose $\kappa$ is surjective, and let $\mathfrak{p} \in \operatorname{Spec}R$.
Since $S$ is a finitely generated $k$-algebra, $R$ is a countably generated $k$-algebra.
By assumption $k$ is uncountable, and thus $R$ is Jacobson.
Therefore, since $\mathfrak{p}$ is prime,
$$\mathfrak{p} = \ \ \ \bigcap_{\mathclap{\mathfrak{p} \subseteq \mathfrak{m} \in \operatorname{Max}R}} \mathfrak{m}.$$
Since $\kappa$ is surjective, the ideal 
$$\mathfrak{q} := \ \ \ \bigcap_{\mathclap{\substack{\mathfrak{n} \in \iota^{-1}(\mathfrak{m}) \text{ s.t.} \\ \mathfrak{p} \subseteq \mathfrak{m} \in \operatorname{Max}R}}} \mathfrak{n}$$
satisfies $\mathfrak{q} \cap R = \mathfrak{p}$.
Furthermore, $\mathfrak{q}$ is radical since it is the intersection of radical ideals.
Thus, since $S$ is noetherian, the Lasker-Noether theorem implies 
$$\mathfrak{q} = \mathfrak{q}_1 \cap \cdots \cap \mathfrak{q}_{\ell},$$ 
where each $\mathfrak{q}_i$ is a minimal prime over $\mathfrak{q}$.

Set $\mathfrak{p}_j := \mathfrak{q}_j \cap R$.
Since $\mathfrak{p} = \mathfrak{q} \cap R \subseteq \mathfrak{q}_j \cap R = \mathfrak{p}_j$, we have $\mathfrak{p} \subseteq \mathfrak{p}_j$.
Thus
$$\mathfrak{p} = \sqrt{ \mathfrak{p}^{\ell}} \subseteq \sqrt{ \mathfrak{p}_1 \cdots \mathfrak{p}_{\ell}}  \subseteq 
\sqrt{ \mathfrak{p}_1 \cap \cdots \cap \mathfrak{p}_{\ell}} = \sqrt{ \left( \mathfrak{q}_1 \cap \cdots \cap \mathfrak{q}_{\ell} \right) \cap R} = \sqrt{ \mathfrak{q} \cap R} = \mathfrak{p}.$$
Whence
$$\mathfrak{p} = \sqrt{\mathfrak{p}_1 \cap \cdots \cap \mathfrak{p}_{\ell}}.$$
Therefore, since $\mathfrak{p}$ is prime and each $\mathfrak{p}_j$ is prime (Lemma \ref{subset}), there is some $1 \leq i \leq \ell$ such that $\mathfrak{p} = \mathfrak{p}_i$. 
Indeed, otherwise each $\mathfrak{p}_i$ would contain some $a_i$ not in $\mathfrak{p}$, but this is not possible since then $a_1 \cdots a_{\ell} \in \mathfrak{p}$ and $\mathfrak{p}$ is prime.\footnote{In general, $\mathfrak{p}$ need not equal $\mathfrak{p}_j$.  Indeed, consider $S = k[x,y]$, $R = k + xS$, and the prime ideals
$$\mathfrak{p} = (xy, xy^2, \ldots)R \in \operatorname{Spec}R, \ \ \ \mathfrak{q}_1 = xS \in \operatorname{Spec}S, \ \ \ \mathfrak{q}_2 = yS \in \operatorname{Spec}S.$$
Then 
$$\mathfrak{p} = \left( \mathfrak{q}_1 \mathfrak{q}_2 \right) \cap R = xyS \cap R = \sqrt{ \left( xS \cap R \right) \left( yS \cap R \right)} = \sqrt{ \left( \mathfrak{q}_1 \cap R \right) \left( \mathfrak{q}_2 \cap R \right)}$$
and $\mathfrak{p}_2 := \mathfrak{q}_2 \cap R = \mathfrak{p}$, but $\mathfrak{p}_1 := \mathfrak{q}_1 \cap R \not = \mathfrak{p}$.} 
Therefore 
$$\mathfrak{q}_i \cap R = \mathfrak{p}_i = \mathfrak{p}.$$
It follows that $\iota$ is surjective.
\end{proof}

\begin{Lemma} \label{ht p} 
Suppose $S$ is a depiction of $R$.
\begin{enumerate}
 \item If $\mathfrak{p} \in \operatorname{Spec}R$ and $\mathfrak{q} \in \iota^{-1}(\mathfrak{p})$, then $\operatorname{ht}\left( \mathfrak{q} \right) \leq \operatorname{ht}\left(\mathfrak{p} \right)$.
 \item If $\mathfrak{m} \in \operatorname{Max}R$, then $\operatorname{ht}(\mathfrak{m}) = \operatorname{dim}R$.
\end{enumerate}
\end{Lemma}

\begin{proof}
(1) Let $\mathfrak{p} \in \operatorname{Spec}R$, $\mathfrak{q} \in \iota^{-1}(\mathfrak{p})$, and set $U := U_{S/R}$.
If $\mathcal{Z}(\mathfrak{q}) \cap U \not = \emptyset$, then $S_{\mathfrak{q}} = R_{\mathfrak{p}}$ by Lemma \ref{Rp Sq}.
Whence 
$$\operatorname{ht}\left( \mathfrak{q} \right) = \operatorname{dim} S_{\mathfrak{q}} = \operatorname{dim} R_{\mathfrak{p}} = \operatorname{ht}(\mathfrak{p}).$$

So suppose $\mathcal{Z}(\mathfrak{q}) \cap U = \emptyset$.
Since $U \not = \emptyset$, $U$ is an open dense set by Proposition \ref{electric}. 
Therefore there is a maximal chain of prime ideals in $S$ containing $\mathfrak{q}$,
$$0 \subsetneq \mathfrak{q}_1 \subsetneq \cdots \subsetneq \mathfrak{q}_{\ell -1} \subsetneq \mathfrak{q}_{\ell} = \mathfrak{q},$$
such that $\mathcal{Z}\left(\mathfrak{q}_i \right) \cap U \not = \emptyset$ for each $1 \leq i \leq \ell -1$.

Set $\mathfrak{p}_i := \mathfrak{q}_i \cap R \in \operatorname{Spec}R$.
Then by Lemma \ref{Rp Sq}, $R_{\mathfrak{p}_i} = S_{\mathfrak{q}_i}$.
In particular, 
$$\operatorname{ht}\left( \mathfrak{p}_{\ell -1} \right) = \operatorname{dim} R_{\mathfrak{p}_{\ell -1}} = \operatorname{dim} S_{\mathfrak{q}_{\ell -1}} = \operatorname{ht}\left( \mathfrak{q}_{\ell -1} \right).$$
Furthermore, by the contrapositive of Lemma \ref{qq'}, 
$$\mathfrak{p}_{\ell -1} \subsetneq \mathfrak{p}.$$
Therefore $\operatorname{ht}\left( \mathfrak{q} \right) \leq \operatorname{ht}\left( \mathfrak{p} \right)$.

(2) Let $\mathfrak{m} \in \operatorname{Max}R$.
Since $\iota$ is surjective, there is some $\mathfrak{q} \in \operatorname{Spec}S$ such that $\iota(\mathfrak{q}) = \mathfrak{m}$.
Let $\mathfrak{n}$ be a maximal ideal of $S$ containing $\mathfrak{q}$.
Then $\mathfrak{m} = \mathfrak{q} \cap R \subseteq \mathfrak{n} \cap R$, and so $\mathfrak{m} = \mathfrak{n} \cap R$ since $\mathfrak{m}$ is maximal.
Therefore
$$\operatorname{dim}R \stackrel{\textsc{(i)}}{=} \operatorname{dim}S \stackrel{\textsc{(ii)}}{=} \operatorname{ht}(\mathfrak{n}) \stackrel{\textsc{(iii)}}{\leq} \operatorname{ht}(\mathfrak{m}) \leq \operatorname{dim}R.$$ 
Indeed, (\textsc{i}) holds by Theorem \ref{U and W}.4; (\textsc{ii}) holds since $S$ an integral domain; and (\textsc{iii}) holds by Claim (1).
\end{proof}

\begin{Theorem} \label{codim theorem}
Suppose $R$ admits a depiction and let $\mathfrak{p} \in \operatorname{Spec}R$.
Then 
\begin{equation} \label{codim circ}
\operatorname{ght}(\mathfrak{p}) \leq \operatorname{ht}(\mathfrak{p}),
\end{equation}
with equality if there is a depiction $S$ of $R$ and $\mathfrak{q} \in \iota^{-1}_S(\mathfrak{p})$ such that $\mathcal{Z}\left( \mathfrak{q} \right) \cap U_{S/R} \not = \emptyset$.
Furthermore,
\begin{enumerate} 
 \item If $R$ is noetherian, then 
$$\operatorname{ght}(\mathfrak{p}) = \operatorname{ht}(\mathfrak{p}).$$ 
 \item If $\mathfrak{m} \in \operatorname{Max}R$, $\mathfrak{q} \in \iota^{-1}_S(\mathfrak{m})$, and $\operatorname{dim} \mathcal{Z}\left( \mathfrak{q} \right) \geq 1$, then 
 $$\operatorname{ght}(\mathfrak{m}) \not = \operatorname{ht}(\mathfrak{m}).$$
\end{enumerate}
\end{Theorem}

\begin{proof}
The inequality (\ref{codim circ}) holds by Lemma \ref{ht p}.1.
Furthermore, if $\mathcal{Z}(\mathfrak{q}) \cap U_{S/R} \not = \emptyset$, then $\operatorname{ght}(\mathfrak{p}) = \operatorname{ht}(\mathfrak{p})$ by Lemma \ref{nonlocal points}.

Suppose $R$ is noetherian.
Then $R$ is a depiction of itself with $U_{R/R} = \operatorname{Max}R$.
Therefore Claim (1) holds as a particular case of the previous paragraph.

Now assume the hypotheses of Claim (2).
Then
$$\operatorname{ght}(\mathfrak{m})
\stackrel{\textsc{(i)}}{\leq}
\operatorname{ht}(\mathfrak{q})
\stackrel{\textsc{(ii)}}{<}
\operatorname{dim}S
\stackrel{\textsc{(iii)}}{=}
\operatorname{dim}R
\stackrel{\textsc{(iv)}}{=}
\operatorname{ht}(\mathfrak{m}).$$
Indeed, (\textsc{i}) holds by Definition \ref{dimension def}; (\textsc{ii}) holds since $\operatorname{dim}\mathcal{Z}(\mathfrak{q}) \geq 1$; (\textsc{iii}) holds by Theorem \ref{U and W}.4; and (\textsc{iv}) holds by Lemma \ref{ht p}.2.
\end{proof}

A depiction, by definition, is a finitely generated $k$-algebra.
However, to capture a similar notion for local rings we introduce the following.

\begin{Definition} \rm{
Suppose $S$ is a depiction of $R$ and let $\mathfrak{q} \in \operatorname{Spec}S$.
Then we say $S_{\mathfrak{q}}$ is a \textit{local depiction} of $R_{\mathfrak{q} \cap R}$.
}\end{Definition}

\begin{Remark} \rm{
A local depiction may not satisfy conditions (i) or (iii) in Definition \ref{dimension def}.
For example, consider $S = k[x,y]$ and $R = k + xS$.
Then the local ring $S_{xS}$ is a noetherian overring of the local ring $R_{xS}$, but $U_{S_{xS}/R_{xS}} = \emptyset$.

Furthermore, although $\iota_S: \operatorname{Spec}S \to \operatorname{Spec}R$ is surjective, the morphism $\iota_{S_{xS}}: \operatorname{Spec}S_{xS} \to \operatorname{Spec}R_{xS}$ is not surjective.
Indeed, there are only two prime ideals of $S$ contained in $xS$, 
\begin{equation} \label{be sure}
0 \subsetneq xS.
\end{equation}
However, there are three prime ideals of $R$ contained in $xS$, namely
\begin{equation} \label{to drink your ovaltine}
0 \subsetneq yS \cap R \subsetneq xS \cap R = xS.
\end{equation}
(These ideals are prime by Lemma \ref{subset}.)
Therefore 
$$\iota_{S_{xS}}: \operatorname{Spec}S_{xS} \to \operatorname{Spec}R_{xS}$$ 
is not surjective. 
}\end{Remark}

\begin{Example} \label{ho ho ho} \rm{
We give an explicit example where the height and geometric codimension of an ideal do not coincide.
Let $S = k[x,y]$ and $R = k + xS$.
Then the chains (\ref{be sure}) and (\ref{to drink your ovaltine}) imply
$$\operatorname{ght} (xS \cap R) = 1 \not = 2 = \operatorname{ht}(xS \cap R).$$ 
}\end{Example}

The following two theorems establish relationships between depictions, geometric dimension, and noetherianity.

\begin{Theorem} \label{non-max}  
Suppose $S$ is a depiction of $R$.
Then the following are equivalent:
 \begin{enumerate} 
  \item $R$ is noetherian.
  \item $U_{S/R} = \operatorname{Max}S$.
  \item $R = S$.
  \end{enumerate}
 In particular, if $R$ is noetherian, then its only depiction is itself.
\end{Theorem}

\begin{proof}
($2 \Rightarrow 1, 3$) 
Suppose $U_{S/R} = \operatorname{Max}S$.
Then
$$R \stackrel{\textsc{(i)}}{=} \bigcap_{\mathfrak{m} \in \operatorname{Max}R} R_{\mathfrak{m}}
\stackrel{\textsc{(ii)}}{=} \bigcap_{\mathfrak{n} \in \operatorname{Max}S} R_{\mathfrak{n} \cap R}
\stackrel{\textsc{(iii)}}{=} \bigcap_{\mathfrak{n} \in \operatorname{Max}S} S_{\mathfrak{n}}
\stackrel{\textsc{(iv)}}{=} S.$$
Indeed, (\textsc{i}) and (\textsc{iv}) hold since $R$ and $S$ are unital commutative rings; (\textsc{ii}) holds since $S$ is a depiction of $R$, whence $\kappa: \operatorname{Max}S \to \operatorname{Max}R$ is surjective; and (\textsc{iii}) holds since $U_{S/R} = \operatorname{Max}S$.
Therefore $R = S$ is noetherian.

($2 \Rightarrow 3$) 
Alternatively, suppose there is some $g \in S \setminus R$.
Then the fractional ideal $(R:g)_R$ is proper, and is thus contained in some maximal ideal $\mathfrak{m} \in \operatorname{Max}R$.
Let $\mathfrak{n} \in \kappa^{-1}(\mathfrak{m})$.
Then $g \in S_{\mathfrak{n}} \setminus R_{\mathfrak{m}}$.
Thus $\mathfrak{n} \in U^c_{S/R}$, and therefore $U^c_{S/R} \not = \emptyset$.

($1 \Rightarrow 2$)
Suppose $R$ is noetherian.
Then $R$ is a depiction of itself.
Let $S$ be another depiction of $R$.
Then
$$\iota_S(U_{S/R}) \subseteq \iota_S\left( \operatorname{Max}S \right) \subseteq \operatorname{Max}R = \iota_R(U_{R/R}) \stackrel{(\textsc{i})}{=} \iota_S(U_{S/R}),$$
where (\textsc{i}) holds by Lemma \ref{nonlocal points}.
Whence $\iota_S(U_{S/R}) = \iota_S \left( \operatorname{Max}S \right)$.
But $\iota_S$ is injective on $U_{S/R}$ by Theorem \ref{U and W}.1.
Therefore $U_{S/R} = \operatorname{Max}S$.

($3 \Rightarrow 1$) If $R = S$, then $R$ is trivially noetherian.
\end{proof}

\begin{Theorem} \label{pos dim theorem}
If $\operatorname{Max}R$ contains a closed point of positive geometric dimension, then $R$ is nonnoetherian.
Conversely, if $R$ is nonnoetherian, depicted by $S$, and there is some $\mathfrak{m} \in \iota(U^c_{S/R})$ such that $\sqrt{\mathfrak{m}S} = \mathfrak{m}$, then $\operatorname{Max}R$ contains a closed point of positive geometric dimension.
\end{Theorem}

\begin{proof}
(i) Suppose $R$ contains a closed point of positive geometric dimension.
Then $R$ admits a depiction $S$.
Let $I$ be a radical ideal of $S$ such that $I \cap R =: \mathfrak{m}$ is a maximal ideal of $R$ and $\operatorname{dim}\mathcal{Z}(I) \geq 1$.
Let $\mathfrak{q} \in \operatorname{Spec}S$ be a minimal prime over $I$.
Then $\operatorname{dim}\mathcal{Z}(\mathfrak{q}) \geq 1$.
Furthermore, $\mathfrak{m} = I \cap R \subseteq \mathfrak{q} \cap R \not = R$ implies $\mathfrak{m} = \mathfrak{q} \cap R$ since $\mathfrak{m}$ is maximal.
Thus by Theorem \ref{codim theorem}.2,
$$\operatorname{ght}(\mathfrak{m}) \not = \operatorname{ht}(\mathfrak{m}).$$
Therefore $R$ is not noetherian by Theorem \ref{codim theorem}.1.

(ii) Suppose the hypotheses hold, and $\operatorname{dim}U^c_{S/R} = 0$.
We claim that $R$ is noetherian.
More specifically, we claim that $R$ is a finitely generated $k$-algebra.
To show this, it suffices to show that $S$ is a finitely generated $R$-module by the Artin-Tate lemma. 

By the Lasker-Noether theorem, there are ideals $\mathfrak{n}_1, \ldots, \mathfrak{n}_{\ell} \in \operatorname{Spec}S$ such that 
$$\mathfrak{m} = \sqrt{\mathfrak{m}S} = \mathfrak{n}_1 \cap \cdots \cap \mathfrak{n}_{\ell}.$$
Whence
$$\mathfrak{m} = \mathfrak{m} \cap R = \left( \mathfrak{n}_1 \cap \cdots \cap \mathfrak{n}_{\ell} \right) \cap R = (\mathfrak{n}_1 \cap R) \cap \cdots \cap (\mathfrak{n}_{\ell} \cap R).$$
Since $\mathfrak{m}$ is maximal in $R$ and each $\mathfrak{n}_i \cap R$ is proper, we have
$$\mathfrak{n}_1 \cap R = \cdots = \mathfrak{n}_{\ell} \cap R = \mathfrak{m}.$$
Therefore, since $\operatorname{dim}U^c_{S/R} = 0$, each $\mathfrak{n}_i$ is a maximal ideal of $S$.

Let $x_1, \ldots, x_t$ be a minimal generating set for $S$ over $k$.
Since $k$ is algebraically closed, for each $1 \leq i \leq \ell$ and $1 \leq j \leq t$ there are scalars $\alpha_{ij} \in k$ such that
\begin{equation} \label{alphaij}
x_j - \alpha_{ij} \in \mathfrak{n}_i.
\end{equation}
Set the degree of each $x_j$ to be 1.
Denote by $\mathcal{M}$ the set of monomials in the variables $x_1, \ldots, x_t$, with coefficient 1 and degree at most $\ell -1$. 
Since $S$ is a finitely generated $k$-algebra, $\mathcal{M}$ is finite.

We proceed by induction to show that as $k$-spaces,
\begin{equation} \label{induction claim}
S = R + \sum_{m \in \mathcal{M}} km.
\end{equation}

\textit{Base case.}
Consider a monomial $x_{j_1} \cdots x_{j_{\ell}}$ of degree $\ell$.
Set 
$$r := \left(x_{j_1} - \alpha_{1j_1} \right)\left(x_{j_2} - \alpha_{2j_2} \right) \cdots \left( x_{j_{\ell}} - \alpha_{\ell j_{\ell}} \right) 
\ \ \ \text{ and } \ \ \
h := x_{j_1}\cdots x_{j_{\ell}} - r.$$ 
Then $r$ is in $\mathfrak{n}_1 \cap \cdots \cap \mathfrak{n}_{\ell} = \mathfrak{m} \subset R$ by (\ref{alphaij}).
Furthermore, $h$ has degree at most $\ell -1$.
Thus $h$ is in $\sum_{m \in \mathcal{M}} k m$ by the definition of $\mathcal{M}$.
Therefore $x_{j_1} \cdots x_{j_{\ell}} = r + h$ is in $R + \sum_{m \in \mathcal{M}}km$.

\textit{Inductive step.}
Now suppose all monomials of degree at most $d-1$ are in $R + \sum_{m \in \mathcal{M}} km$, and consider a monomial $x_{j_1} \cdots x_{j_d}$ of degree $d$.
Set
$$r := \left(x_{j_1} - \alpha_{1j_1} \right) \cdots \left( x_{j_{\ell}} - \alpha_{\ell j_{\ell}} \right)x_{j_{\ell+1}} \cdots x_{j_d} \ \ \ \text{ and } \ \ \ h := x_{j_1}\cdots x_{j_{d}} - r.$$
Then $r$ is in $\mathfrak{n}_1 \cap \cdots \cap \mathfrak{n}_{\ell} = \mathfrak{m} \subset R$, again by (\ref{alphaij}). 
Furthermore, $h$ has degree at most $d-1$.
Thus $h$ is in $R + \sum_{m \in \mathcal{M}} km$ by the induction hypothesis.
Therefore $x_{j_1} \cdots x_{j_d} = r + h$ is in $R + \sum_{m \in \mathcal{M}}km$.
This proves our claim (\ref{induction claim}).

But 
$$R + \sum_{m \in \mathcal{M}} km \subseteq R + \sum_{m \in \mathcal{M}}Rm.$$
Thus as $R$-modules,
$$S = R + \sum_{m \in \mathcal{M}}Rm.$$
Therefore $S$ is a finitely generated $R$-module since $\mathcal{M}$ is finite.
\end{proof}

The following corollary is immediate.

\begin{Corollary} \label{R = k + I}
Let $I$ be a radical ideal of a finitely generated $k$-algebra $S$. 
Then the ring $R = k + I$ is nonnoetherian if and only if $\operatorname{dim}\mathcal{Z}(I) \geq 1$.
\end{Corollary}

\begin{Example} \rm{
Let $S$ be a finitely generated $k$-algebra, and let $\mathfrak{n}_1, \mathfrak{n}_2, \ldots, \mathfrak{n}_{\ell}$ be a finite set of maximal ideals of $S$.
Then by Corollary \ref{R = k + I}, the ring
$$R = k + \sqrt{ \mathfrak{n}_1 \mathfrak{n}_2 \cdots \mathfrak{n}_{\ell}}$$
is noetherian.
Furthermore, the $\ell$ points $\mathfrak{n}_1, \ldots, \mathfrak{n}_{\ell}$ of $\operatorname{Max}S$ are identified as one point in $\operatorname{Max}R$, and all other points are unaltered by Proposition \ref{I},
$$U_{S/R} = \mathcal{Z}(\mathfrak{n}_1 \cdots \mathfrak{n}_{\ell})^c.$$
In particular, the variety $\operatorname{Max}R$ is nonlocal.
$R$ is therefore an example of a noetherian ring with nonlocal geometry, although its one nonlocal point is zero-dimensional.
}\end{Example}

The following proposition characterizes the assumptions in the converse implication of Theorem \ref{pos dim theorem}. 

\begin{Proposition} \label{isolated sing}
If there is a point $\mathfrak{m} \in \iota(U^c_{S/R})$ satisfying $\mathfrak{m}S = \mathfrak{m}$, then
\begin{equation} \label{Uc =}
U^c_{S/R} = \left\{ \mathfrak{m} \right\}.
\end{equation}
If additionally $S$ is a depiction of $R$, then $R$ is an isolated nonnoetherian singularity.
\end{Proposition}

\begin{proof}
Suppose $\mathfrak{m} = \mathfrak{m}S$.
Consider $\mathfrak{m}' \in \operatorname{Max}R \setminus \left\{ \mathfrak{m} \right\}$.
Fix $\mathfrak{n}' \in \iota^{-1}(\mathfrak{m}')$ and $g \in \mathfrak{m} \setminus \mathfrak{m}'$.
Since $\mathfrak{m} = \mathfrak{m}S$, the ideal $Sg$ is in $\mathfrak{m} \subset R$.
Thus
$$S = Sg \cdot g^{-1} \subset R_{\mathfrak{m}'}.$$
Whence $R_{\mathfrak{m}'} = S_{\mathfrak{n}'}$.
Thus $\mathfrak{m}' \in \iota(U_{S/R})$.
Therefore (\ref{Uc =}) holds. 

If $S$ is a depiction of $R$, then $U^*_{S/R} = U_{S/R}$, and so (\ref{Uc =}) implies that $R$ is an isolated nonnoetherian singularity.
\end{proof}

\begin{Question} \rm{
Is there a nonnoetherian algebra (which admits a depiction) such that all of its closed points have geometric dimension zero? 
By Proposition \ref{isolated sing}, such an algebra would necessarily be a non-isolated nonnoetherian singularity.
Or is it the case that $R$ is nonnoetherian if and only if $\operatorname{Max}R$ contains a closed point of positive geometric dimension?
}\end{Question}

\begin{Example} \rm{
Consider the algebras
$$S = k\left[x,y \right], \ \ \ R = k + xS, \ \ \ R' = k\left[x,xy \right].$$
Set $\mathfrak{n} := (x,y)S$.
By Example \ref{ho ho ho}, the closed point $\mathfrak{n} \cap R$ of $\operatorname{Max}R$ has geometric dimension 1.
Naively it appears that the closed point $\mathfrak{n} \cap R'$ of $\operatorname{Max}R'$ should also have geometric dimension 1, contrary to the claim of Theorem \ref{pos dim theorem} since $R'$ is noetherian.
However, although $R$ is depicted by $S$, $R'$ is not.
Indeed, the morphism $\iota: \operatorname{Spec}S \to \operatorname{Spec}R'$ is not surjective: the ideal $xR'$ is prime in $R'$, but 
$$\iota^{-1}(xR') = \emptyset.$$
In contrast, the ideal $xR$ is \textit{not} prime in $R$ since $x \cdot xy^2 = xy \cdot xy$.

In view of Lemma \ref{surjective1}, consider the maximal ideals 
$$\mathfrak{m}_{\alpha} := (x,xy-\alpha)R' \in \operatorname{Max}R'$$ 
with $\alpha \in k$.
Then 
$$xR' = \cap_{\alpha \in k} \mathfrak{m}_{\alpha}.$$
Furthermore, if $\alpha \not = 0$, then $1 \in \mathfrak{m}_{\alpha}S$.
Thus if $\mathfrak{n} \in \operatorname{Max}S$ satisfies $\mathfrak{n} \cap R = \mathfrak{m}_{\alpha}$, then 
$$1 \in \mathfrak{m}_{\alpha}S \subseteq \mathfrak{n},$$ 
which is not possible.
Therefore the morphism $\kappa: \operatorname{Max}S \to \operatorname{Max}R'$ is also not surjective.
}\end{Example}

We say a depiction $S$ of $R$ is \textit{maximal} (resp.\ \textit{minimal}) if $S$ is not contained in (resp.\ does not contain) any other depiction of $R$. 

\begin{Proposition} \label{min and max} \ 
\begin{enumerate}
 \item Minimal depictions do not exist in general.
 \item Maximal depictions are not unique in general.
\end{enumerate}
\end{Proposition}

\begin{proof}
(1) We first show that minimal depictions need not exist.
Let $S = k[x,y]$ and $R = k + xS$.
For $\ell \in \mathbb{N}$, set 
$$S_{\ell} := R\left[y^{\ell}, y^{\ell+1}, y^{\ell+2}, \ldots \right] = k\left[x,xy,xy^2, \ldots, xy^{\ell-1}, y^{\ell}, y^{\ell+1}, \ldots, y^{2\ell-1} \right].$$
Then each $S_{\ell}$ is a depiction of $R$.
However,
$$S_{\ell+1} \subsetneq S_{\ell} \ \ \ \text{ and } \ \ \ R = \bigcap_{\ell \geq 1} S_{\ell}.$$

(2) We now show that maximal depictions need not be unique. 
Let $T = k[x,y,z]$ and $R = k + xyT$.
We claim that the overrings 
$$S := T[x^{-1}] \ \ \text{ and } \ \ S' := T[y^{-1}]$$
are both depictions of $R$.

Indeed, $U_{S/R}$ is nonempty: 
Let $a,b \in k^*$.
Then $x^{-1}$ and each $f \in T$ are in $R_{(x-a,y-b,z)S \cap R}$ since
$$x^{-1} = xy^2 \cdot \frac{1}{(xy)^2} \ \ \text{ and } \ \ f = xyf \cdot \frac{1}{xy}.$$
Thus
$$R_{(x-a,y-b,z)S \cap R} = S_{(x-a,y-b,z)S}.$$
Therefore the maximal ideal $(x-a,y-b,z)S$ is in $U_{S/R}$.
Similarly $U_{S'/R}$ is nonempty.

It is straightforward to check that $\iota_S: \operatorname{Spec}S \to \operatorname{Spec}R$ and $\iota_{S'}: \operatorname{Spec}S' \to \operatorname{Spec}R$ are both surjective by Lemma \ref{surjective1} and noting that 
$$yS \cap R = xS' \cap R = xyT \in \operatorname{Max}R.$$

Finally, the minimal proper overring $T[x^{-1}, y^{-1}]$ of $S$ and $S'$ is not a depiction of $R$: $T[x^{-1}, y^{-1}]$ has no maximal ideal $\mathfrak{n}$ satisfying $\mathfrak{n} \cap R = xyT$.
Therefore 
$$\iota_{T\left[ x^{-1},y^{-1} \right]}: \operatorname{Spec}T\left[x^{-1},y^{-1} \right] \to \operatorname{Spec}R$$
is not surjective.
\end{proof}

\section{Noncommutative nonnoetherian geometry} \label{depictions from impressions}

Throughout $k$ is an uncountable algebraically closed field of characteristic zero; $A = kQ/I$ is a quiver algebra with finite quiver $Q$ and $I \subset kQ_{\geq 1}$; and $B$ is a finitely generated integral domain over $k$.
We say an element $c \in A = kQ/I$ is a cycle (resp.\ path) if there is a cycle (resp.\ path) $c' \in kQ$ such that $c = c' + I$.

Denote by $E_{ji} \in M_{|Q_0|}(B)$ the matrix with a 1 in the $(ji)$-th slot and zeros elsewhere. 
Let $\tau: A \rightarrow M_{|Q_0|}(B)$ be an algebra homomorphism such that $\tau(e_i) = E_{ii}$ for each $i \in Q_0$.
For $p \in e_jAe_i$, denote by $\bar{\tau}(p) \in B$ the single non-zero entry of the matrix $\tau(p)$.
For each $i,j \in Q_0$, $\bar{\tau}$ defines a $k$-linear map $\bar{\tau}: e_jAe_i \rightarrow B$. 
Set
$$\begin{array}{c}
R := k\left[ \displaystyle{\cap_{i \in Q_0}} \bar{\tau}\left( e_iAe_i \right) \right] \subseteq B,\\
S := k \left[ \cup_{i \in Q_0} \bar{\tau}\left( e_iAe_i \right) \right] \subseteq B.
\end{array}$$

Recall the definition of an impression given in Definition \ref{impression definition}. 
An impression is useful in part because it determines the center $Z$ of $A$ \cite[Lemma 2.1]{B}:
\begin{equation} \label{ZR}
Z \cong \left\{ f \in B \ | \ f1_d \in \operatorname{im}\tau \right\} \subseteq B.
\end{equation}
  
\begin{Theorem} \label{impression corollary}
Suppose $\tau: A \hookrightarrow M_{|Q_0|}(B)$ is an impression of $A$ with $B$ an integral domain and $\tau(e_i) = E_{ii}$ for each $i \in Q_0$.
Then $U_{S/R} \not = \emptyset$.
Furthermore, if $U^*_{S/R} \subseteq U_{S/R}$, then
\begin{enumerate}
 \item The center $Z$ of $A$ is isomorphic to $R$ and is depicted by $S$.
 \item The statements
  \begin{enumerate}
   \item $R = S$.
   \item $A$ is a finitely generated $Z$-module.
   \item $Z$ is noetherian.
   \item $A$ is noetherian.
  \end{enumerate}
  satisfy the equivalences
  $$\xy 
 (-24,0)*+{\text{(a)}}="1";(-8,0)*+{\text{(b)}}="2";(8,0)*+{\text{(c)}}="3";(24,0)*+{\text{(d)}}="4";
{\ar@{<=>}"1";"2"};{\ar@{<=>}"2";"3"};{\ar@{=>}@/^/"3";"4"};{\ar@{=>}@/^/^{\star}"4";"3"};
  \endxy$$
  where $(\star)$ holds if the $\bar{\tau}$-image of a path is a monomial in $B$ and $I$ is generated by binomials in the paths of $Q$.
\end{enumerate}
\end{Theorem}

In Example \ref{y2 - ba} below we show that the assumption $U^*_{S/R} \subseteq U_{S/R}$ is independent from the impression assumption.

\begin{proof}
(1.i) We first show that the center of $A$ is isomorphic to $R$.
Set $d := |Q_0|$.
Clearly 
$$\left\{ f \in B \ | \ f1_d \in \operatorname{im}\tau \right\} \subseteq R.$$
So suppose $f \in R$.
Then for each $i \in Q_0$ there is some $c_i \in e_iAe_i$ such that $\bar{\tau}(c_i) = f$.
Thus $\tau\left( \sum_{i \in Q_0} c_i \right) = f1_d$.
Whence $\left\{ f \in B \ | \ f1_d \in \operatorname{im}\tau \right\} \supseteq R$.
Thus by (\ref{ZR}),
$$Z \cong \left\{ f \in B \ | \ f1_d \in \operatorname{im}\tau \right\} = R.$$

(1.ii) We now show that $R$ is depicted by $S$. 

(1.ii.\textsc{i}) \textit{$S$ is a finitely generated domain over $k$.} 
By assumption $Q$ is finite.
Thus there are a finite number of non-vertex cycles in $Q$ without cyclic proper subpaths (not viewed modulo $I$), say $c_1, \ldots, c_{\ell}$.
Note that each cycle has length at most $|Q_0|$.

We claim that 
\begin{equation} \label{S = kc}
S = k\left[ \bar{\tau}\left(c_1\right), \ldots, \bar{\tau}\left(c_{\ell} \right) \right].
\end{equation}
It suffices to show that the $\bar{\tau}$-image of each non-vertex cycle with a cyclic proper subpath is contained in $k\left[ \bar{\tau}\left(c_1\right), \ldots, \bar{\tau}\left(c_{\ell} \right) \right].$
So suppose $d$ is such a cycle.
Then $d$ has a cyclic proper subpath $c$ with no cyclic proper subpaths (again, not viewed modulo $I$).
Thus there are paths $d_1$, $d_2$ such that $d = d_2cd_1$.
Since $c$ is a cycle, $d_2d_1$ is also a cycle.
Therefore, since $\tau$ is an algebra homomorphism and $B$ is commutative,
$$\bar{\tau}(d) = \bar{\tau}\left(d_2cd_1\right) = \bar{\tau}\left(d_2\right) \bar{\tau}\left(c \right) \bar{\tau}\left(d_1\right) = \bar{\tau}\left(c\right) \bar{\tau}\left(d_2\right)\bar{\tau}\left(d_1\right) = \bar{\tau}\left(c\right) \bar{\tau}\left(d_2d_1\right).$$
The length of the cycle $d_2d_1$ is strictly less than the length of $d$ since $c$ is a non-vertex cycle.
Our claim (\ref{S = kc}) then follows by induction on the length of the cycles.

Furthermore, $S$ is a domain since it is a subalgebra of the domain $B$.

(1.ii.\textsc{ii}) 
\textit{The set $U_{S/R}$ is nonempty.} 
By \cite[Lemma 2.4]{B}, the dimension vector for the simple $A$-modules of maximal $k$-dimension is $1^{Q_0}$.
Thus there exists a path $p_{ji} \not \in I$ between any two vertices $i,j$ of $Q$.
The cycle $b := p_{1 d} \cdots p_{32}p_{21}$ then contains each vertex as a subpath.
Furthermore, since $\tau$ is injective, $\bar{\tau}(p_{ji}) \not = 0$.
Thus, since $B$ is an integral domain and $\tau$ is an algebra homomorphism,
$$\bar{\tau}(b) = \bar{\tau}\left( p_{1 d} \cdots p_{32}p_{21} \right) = \bar{\tau}\left( p_{1 d} \right) \cdots \bar{\tau}\left( p_{32} \right) \bar{\tau}\left( p_{21} \right) \not = 0.$$
Whence $\bar{\tau}(b) \not = 0$.

Fix $i \in Q_0$, and let $c_i \in e_iAe_i$ be an arbitrary cycle.  
For each $j \in Q_0$, denote by $b_j$ and $d_j$ the respective cycles obtained by cyclically permuting $b$ and $d_i := bc_i$ so that their heads and tails are at $j$.  
Then $\bar{\tau}(b_j) = \bar{\tau}(b) =: \beta$ and $\bar{\tau}(d_j) = \bar{\tau}(d_i) = \bar{\tau}(c_i) \beta$, since $\tau$ is an algebra homomorphism on $e_iAe_i$.
Therefore $\beta$ and $\bar{\tau}(c_i)\beta$ are in $R$.  

Fix $\mathfrak{b} \in \mathcal{Z}(\beta)^c \subset \operatorname{Max}B$.
By Lemma \ref{subset}, $\mathfrak{n}:= \mathfrak{b} \cap S$ and $\mathfrak{m} := \mathfrak{n} \cap R$ are maximal ideals of $S$ and $R$ respectively.
Furthermore, $\beta \in R$ is invertible in the localization $R_{\mathfrak{m}}$.  
Consequently 
$$\bar{\tau}(c_i) = \bar{\tau}(c_i)\beta \cdot \beta^{-1} \in R_{\mathfrak{m}}.$$  
Since $c_i$ was arbitrary, $S_{\mathfrak{n}} \subseteq R_{\mathfrak{m}}$.
Whence $S_{\mathfrak{n}} = R_{\mathfrak{m}}$.   

(1.ii.\textsc{iii})
\textit{The map $\iota: \operatorname{Spec}S \to \operatorname{Spec}R$, $\mathfrak{q} \mapsto \mathfrak{q} \cap R$, is surjective.} 
By Claim (1.ii.a), $S$ is a finitely generated $k$-algebra.
Thus by Lemma \ref{surjective1}, it suffices to show $\kappa: \operatorname{Max}S \to \operatorname{Max}R$, $\mathfrak{n} \mapsto \mathfrak{n} \cap R$, is surjective.
Let $\mathfrak{m} \in \operatorname{Max}R$.  
By the definition of impression, the morphism $\operatorname{Max}B \rightarrow \operatorname{Max}\tau(Z)$, $\mathfrak{b} \mapsto \mathfrak{b}1_d \cap \tau(Z)$, is surjective. 
Thus there is some $\mathfrak{b} \in \operatorname{Max}B$ such that $\mathfrak{b} \cap R = \mathfrak{m}$.
By Lemma \ref{subset}, $\mathfrak{n} := \mathfrak{b} \cap S$ is in $\operatorname{Max}S$.
Furthermore, $\mathfrak{n}$ satisfies 
$$\mathfrak{n} \cap R = (\mathfrak{b} \cap S) \cap R = \mathfrak{m}.$$
Thus $\kappa$ is surjective.

(1.ii.\textsc{iv}) Claims (1.ii.\textsc{i}, \textsc{ii}, \textsc{iii}), together with the assumption $U^*_{S/R} \subseteq U_{S/R}$, imply that $S$ is a depiction of $R$.

(2. a $\Leftrightarrow$ c) Follows from Claim (1) and Theorem \ref{non-max}.2.

(2. a $\Rightarrow$ b,d) Follows from \cite[Theorem 2.11]{B}.

(2. b $\Rightarrow$ a) Suppose $R \not = S$.
Then $R$ is an infinitely generated $k$-algebra by Theorem \ref{non-max}.1.
Furthermore, $S$ is a finitely generated $k$-algebra by Claim (1.ii.\textsc{i}).
Thus $S$ is an infinitely generated $R$-module by the Artin-Tate lemma.
But $\tau$ is injective, $|Q_0| < \infty$, and $Z \cong R$ by Claim (1).
Therefore $\bigoplus_{i \in Q_0} e_iAe_i$ is an infinitely generated $Z$-module.
Whence $A$ is an infinitely generated $Z$-module.

(2. d $\Rightarrow$ a) Suppose $R \not = S$, and the conditions ($\star$) hold.
As was shown in (b $\Rightarrow$ a), $S$ is an infinitely generated $R$-module.
Thus there is a cycle $q$ and vertex $i$ such that for each $n \geq 1$,
$$\bar{\tau}(q)^n \in S \setminus \bar{\tau}\left( e_iAe_i \right).$$

Since $\tau$ is an impression of $A$, there is a path $p_1$ in $e_{\operatorname{t}(q)}Ae_i$ and a path $p_2$ in $e_iAe_{\operatorname{t}(q)}$.
Assume to the contrary that the chain of ideals
$$0 \subsetneq \left(p_2p_1 \right) \subsetneq \left( p_2p_1, p_2qp_1 \right) \subsetneq \left( p_2p_1, p_2qp_1, p_2q^2p_1 \right) \subsetneq \cdots$$
in $e_iAe_i$ terminates.
Then, since $I$ is generated by binomials in the paths of $Q$, there is some $n > m \geq 0$ and a cycle $a \in e_iAe_i$ such that
$$p_2q^np_1 = ap_2q^mp_1.$$
Whence
$$\bar{\tau}(q)^n \bar{\tau}\left(p_2p_1 \right) = \bar{\tau}\left(p_2q^np_1\right) = \bar{\tau}\left(a p_2q^mp_1 \right) = \bar{\tau}(a) \bar{\tau}(q)^m \bar{\tau}\left( p_2p_1\right).$$
Thus since $B$ is an integral domain,
$$\bar{\tau}(q)^{n-m} = \bar{\tau}(a) \in \bar{\tau}\left(e_iAe_i \right).$$
But this contradicts our choice of $q$.

Thus the vertex corner ring $e_iAe_i$ is nonnoetherian.
Therefore $A$ is nonnoetherian.
\end{proof}

\begin{Remark} \rm{
The role of $S$ is new: $S$ is a commutative ring obtained from $A$ that in most cases is not a central subring of $A$, but is closely related to the geometry of the center $Z$ of $A$.  
By Theorem \ref{impression corollary}.2, if $Z$ is noetherian, then $S$ is isomorphic to $Z$, and if $Z$ is nonnoetherian, then $S$ properly contains $Z$.
} \end{Remark}

\begin{Example} \label{warp drive} \rm{
Consider the quiver algebra 
\begin{equation} \label{A = kQ/I}
A = kQ/\left\langle yba - bay \right\rangle
\end{equation}
with quiver given in Figure \ref{first figure}.
$A$ admits the impression $\left( \tau, B = k[x,y] \right)$ where $\tau(e_i) = E_{ii}$ for $i = 1,2$, and
\begin{equation} \label{labels}
\bar{\tau}(a) = 1, \ \ \bar{\tau}(b) = x, \ \ \bar{\tau}(y) = y.
\end{equation}
Recall our motivating example in Section \ref{Section 1}, $S' = k[x,y]$ and $R' = k + xS'$. 
By (\ref{ZR}), the center of $A$ is isomorphic to
$$R = k\left[ \bar{\tau}\left( e_1Ae_1 \right) \cap \bar{\tau}\left(e_2Ae_2 \right) \right] = k + xS' = R'$$
and is depicted by
$$S = k\left[ \bar{\tau}\left( e_1Ae_1 \right) \cup \bar{\tau}\left( e_2Ae_2 \right) \right] = k[x,y] = S'.$$
By Theorem \ref{impression corollary}.2, $A$ and $Z$ are nonnoetherian, and $A$ is an infinitely generated $Z$-module. 
}\end{Example}

\begin{figure}
$$\xy
(-8,0)*+{\text{\scriptsize{$2$}}}="1";(8,0)*+{\text{\scriptsize{$1$}}}="2";
{\ar@/^/^a"1";"2"};{\ar@/^/^b"2";"1"};
{\ar@(ul,dl)_y"1"};
\endxy$$
\caption{}
\label{first figure}
\end{figure}

We recall a homological characterization of smoothness, in the both the commutative and noncommutative settings.
Let $R$ be a noetherian integral domain and let $\mathfrak{p} \in \operatorname{Spec}R$.
Then by the Auslander-Buchsbaum formula,
$$\operatorname{pd}_{R_{\mathfrak{p}}}\left( R_{\mathfrak{p}}/\mathfrak{p} \right) = \left\{ \begin{array}{cl} \operatorname{ht}(\mathfrak{p}) & \text{ if $\mathfrak{p}$ is smooth} \\ \infty & \text{ otherwise}. \end{array} \right.$$
This notion was generalized to the noncommutative setting by Brown and Hajarnavis \cite{BH}.
They define a noetherian (noncommutative) algebra $A$ with prime center $R$ to be \textit{homologically homogeneous} if the projective dimension of each simple $A$-module $V$ equals the Krull dimension of $R$,\footnote{Specifically, if $R$ is a commutative noetherian equidimensional $k$-algebra and $A$ is a module-finite $R$-algebra, then $A$ is homologically homogeneous if all simple $A$-modules have the same projective dimension.} 
\begin{equation} \label{pd ht}
\operatorname{pd}_A(V) = \operatorname{dim}R = \operatorname{ht}(\operatorname{ann}_R(V)).
\end{equation}
Using this notion, Van den Bergh defines a noncommutative crepant resolution $A$ of a noetherian normal Gorenstein domain $R$ to be a homologically homogeneous endomorphism ring $A = \operatorname{End}_R(M)$ of a reflexive finitely generated $R$-module $M$ \cite[Definition 4.1]{V}. 

We propose that if $R$ is a nonnoetherian domain, then homological homogeneity should be replaced by an equality between projective dimension and geometric height, rather than height as in (\ref{pd ht}).
This proposal is illustrated in the following proposition, and is further studied in \cite{B4} in the context of homotopy dimer algebras.

\begin{Proposition} \label{A quiver}
Let $A = kQ/I$ be the quiver algebra (\ref{A = kQ/I}), and set $\mathfrak{m} = xS$.
Then 
\begin{itemize} 
 \item $A$ is an endomorphism of a reflexive module over its center $Z \cong R$,
$$A \cong \operatorname{End}_R\left( R \oplus \mathfrak{m} \right) \cong \operatorname{End}_Z\left(Ae_1\right).$$
 \item Let $V_i$ be the simple $A$-module supported at vertex $i$. 
 Then 
$$\operatorname{pd}_A(V_i) = \operatorname{ght}\left( \operatorname{ann}_{e_iAe_i}(V_i) \right).$$
 \item The smooth locus of $Z$ parameterizes the simple $A$-module isoclasses of maximal $k$-dimension, and coincides with the open set $U_{S/R} \subset \operatorname{Max}S$.
\end{itemize}
Therefore, although $A$ and $Z$ are nonnoetherian and $A$ is an infinitely generated $Z$-module, $A$ nevertheless may be viewed as a noncommutative desingularization of its center.
\end{Proposition}

\begin{proof}
For the following, denote $\operatorname{Hom}_R(-,-)$ by $(-,-)$.
Using the labeling of arrows given by the impression (\ref{labels}), we find
$$\begin{array}{rcl}
\operatorname{End}_Z\left(Ae_1 \right) & = & \operatorname{End}_Z\left( e_1Ae_1 \oplus e_2Ae_1 \right)\\
& \cong & \operatorname{End}_R\left(R \oplus \mathfrak{m} \right)\\
& \cong & \left( \begin{matrix} (R,R) &  (\mathfrak{m},R) \\ (R, \mathfrak{m}) & (\mathfrak{m}, \mathfrak{m}) \end{matrix} \right)\\
& \cong & \left( \begin{matrix} R & S \\ \mathfrak{m} & S \end{matrix} \right)\\
& \cong & \left( \begin{matrix} e_1Ae_1 & e_1Ae_2 \\ e_2Ae_1 & e_2Ae_2 \end{matrix} \right)\\
& \cong & A.
\end{array}$$
Furthermore, $R \oplus \mathfrak{m}$ is a reflexive $R$-module:
$$\left((R \oplus \mathfrak{m}, R), R \right) \cong \left( R \oplus S, R \right) \cong R \oplus \mathfrak{m}.$$

The minimal projective resolution of $V_1$ is
$$0 \rightarrow \left( \begin{matrix} S \\ S \end{matrix} \right)
\stackrel{\cdot x}{\longrightarrow} 
\left( \begin{matrix} R \\ \mathfrak{m} \end{matrix} \right) 
\rightarrow 
\left( \begin{matrix} R/\mathfrak{m} \\ 0 \end{matrix} \right) \cong \left( \begin{matrix} k \\ 0 \end{matrix} \right) = V_1 \rightarrow 0.$$
Set $\mathfrak{n} := (x,y)S \in \operatorname{Max}S$.
Then the minimal projective resolution of $V_2$ is
$$0 \rightarrow \left( \begin{matrix} S \\ S \end{matrix} \right)
\stackrel{\cdot \left( \begin{matrix} xy & -x \end{matrix} \right)}{\longrightarrow} 
\left( \begin{matrix} R & S \\ \mathfrak{m} & S \end{matrix} \right)
\stackrel{\cdot \left( \begin{matrix} 1 \\ y \end{matrix} \right)}{\longrightarrow}
\left( \begin{matrix} S \\ S \end{matrix} \right) 
\rightarrow \left( \begin{matrix} 0 \\ S/\mathfrak{n} \end{matrix} \right) \cong \left( \begin{matrix} 0 \\ k \end{matrix} \right) = V_2 \rightarrow 0.$$

The simple $A$-modules of maximal $k$-dimension are the simples modules with dimension vector $(1,1)$ by \cite[Lemma 5.1]{B}. 
These modules are parameterized by the smooth locus of $\operatorname{Max}Z$, namely $(ab,y) \in k^* \times k$, which coincides with $U_{S/R} \subset \operatorname{Max}S$.
\end{proof}

Although $A$ is isomorphic to $\operatorname{End}_Z(Ae_1)$, note that $A$ is not isomorphic to $\operatorname{End}_Z\left( Ae_2 \right)$ since
$$\operatorname{End}_Z(Ae_2) \cong \operatorname{End}_R(S \oplus S) \cong M_2(S).$$
Furthermore, the moduli space of $\theta$-stable $A$-modules of dimension vector $(1,1)$, for generic stability parameter $\theta$, is precisely the desingularization $\operatorname{Max}S$. 
($\operatorname{Max}S$ is not a resolution of $\operatorname{Max}R$ since the morphism $\kappa: \operatorname{Max}S \to \operatorname{Max}R$ is not proper.)

The following example demonstrates the necessity of the assumption $U^*_{S/R} \subseteq U_{S/R}$ in Theorem \ref{impression corollary}.

\begin{Example} \label{y2 - ba} \rm{
Consider the quiver algebra 
$$A = kQ/\left\langle yba - bay, y^2 - ba \right\rangle$$
with quiver given in Figure \ref{first figure}, as in Example \ref{warp drive}.
$A$ admits an impression $\left(\tau, B = k[x] \right)$ where $\tau(e_i) = E_{ii}$ for $i = 1,2$, and
$$\bar{\tau}(a) = \bar{\tau}(b) = \bar{\tau}(y) = x.$$
By (\ref{ZR}), the center $Z$ of $A$ is isomorphic to 
$$R = k\left[ x^2, x^3\right].$$
Therefore $Z$ is noetherian.
But
$$S = k\left[ x \right] \not = R.$$
However, Theorem \ref{impression corollary}.2 is not applicable to this example because 
$$U^*_{S/R} = \mathbb{A}^1 \not \subseteq  \mathbb{A}^1 \setminus \left\{ 0 \right\} = U_{S/R};$$
see Remark \ref{U not = W}.
}\end{Example}

\ \\
\textbf{Acknowledgments.}  The author would like to thank Paul Smith and Travis Kopp for useful discussions, and an anonymous referee for useful comments. 
Part of this paper is based on work supported by the Simons Foundation while the author was a postdoc at the Simons Center for Geometry and Physics at Stony Brook University. 
His work was also supported in part by a PFGW Grant, which he gratefully acknowledges. 

\bibliographystyle{hep}
\def\cprime{$'$} \def\cprime{$'$}

\end{document}